**Modified convex hull pricing for fixed load power markets**

Vadim Borokhov

LLC "En+development",
ul. Vasilisi Kozhinoy 1, Moscow, 121096, Russia
E-mail: vadimbor@yahoo.com

**Abstract**

We consider a fixed load power market with non-convexities originating from the start-up and no-load costs of generators. The convex hull (minimum-uplift) pricing method results in power prices that minimize the total uplift payment to generators introduced to compensate their potential profits lost by accepting the centralized dispatch solution. In this method, an opportunity to supply any other output volume allowed by the generator private constraints is treated as foregone. For each generator, we identify the output volumes that are attainable at the power market (i.e. the volumes that are both economically and technologically feasible) and define the corresponding modified private feasible set. These sets are further utilized in the lost profit calculations for the generators. The new pricing method results in a generally different set of market prices and lower (or equal) total uplift payment compared to the convex hull pricing algorithm.

## 1. Introduction

Liberalization of the power sector paved the way for development of electricity markets with free pricing for power, which can be either centrally coordinated or decentralized, e.g. based on bilateral trade. The centrally coordinated electricity markets are often based on security-constrained commitment and economic dispatch optimization problems, which are reduced to the least cost commitment and dispatch problems for power systems with the fixed load. If the centralized dispatch optimization problem is convex, then there exists an equilibrium price (not necessarily unique) that supports the solution: given the price, no market player (acting as a price-taker) has economic incentives to distort its output/consumption volumes. The equilibrium price, if it exists, can be obtained by means of the Walrasian auction, where each market player submits its supply (demand) bids. In this framework, the equilibrium price is the price at which the total demand equals the total supply. Given the producer (consumer) cost (benefit) functions, the total supply (demand) volume at a given price is determined using the decentralized dispatch problem obtained by the Lagrangian relaxation of the power balance constraint in the centralized dispatch problem. If the centralized dispatch problem is convex, then an equilibrium price coincides with the marginal price and no producer (consumer) has economic incentives to change its output (consumption) volumes set by the centralized dispatch (Bohn et al. 1984; Littlechild, 1988; Schweppe et al. 1988). In a uninode single-period power market with the convex centralized dispatch problem, the marginal price (the set of marginal prices) is given by the intersection of the aggregate supply curve with the aggregate demand curve. However, if the centralized dispatch problem is non-convex, then the marginal price does not reflect the non-convex features of power output such as nonzero minimum output limits, fixed costs (start-up and no-load costs), as well as other sources of non-convexities. As a result, the marginal price may not compensate the full cost of power output and, therefore, may fail to act as an equilibrium market price. Moreover, the equilibrium price may not even exist and other mechanisms are needed to ensure the stability of the centralized dispatch solution. Many different pricing schemes have been proposed for markets with non-convexities, including the introduction of new products/services (and the associated prices), utilization of a nonlinear pricing methodology (which includes the case of discriminatory, i.e. market player specific, prices and the case of generator/consumer revenues/costs being some nonlinear functions of power volumes), and application of the uniform (linear) pricing for power with the applicable uplifts (side-payments). The new service (a unit being online) and the associated market player specific prices were introduced by O'Neill et al. (2005) by fixing the integral status variables at their optimal values obtained from the solution to the centralized dispatch. However, the pricing scheme results in zero profit for all generators and may produce negative prices for being online. Therefore, such a pricing algorithm is close to pay-as-bid pricing. Moreover, if the negative prices are discarded so that generators retain their profits, then the resulting pricing method does not result in the competitive equilibrium. Bjørndal and Jörnsten (2009) and Bjørndal and Jörnsten (2010) amended this approach to produce more stable prices by adding extra constraints to the reformulated optimization problem that fix certain continuous variables at their optimal values as well. Andrianesis et. al. (2013a), (2013b) presented the compensation mechanisms that cover the losses of the generating units participating in the centralized market and analyzed the resulting bidding incentives for the units.

The discriminatory pricing methods were proposed by Madrigal and Quintana (2001), Madrigal, Quintana, and Aguado (2001), Hao and Zhuang (2003), Santiago Lopez and Madrigal (2005). Motto and Galiana (2002), Galiana, Motto, and Bouffard (2003), Bouffard, and Galiana (2005) introduced the nonlinear pricing in the form of the generalized uplift functions that include generators and consumers in the lost profit compensation and result in zero-sum uplift at the market. Liberopoulos and Andrianesis (2016) proposed the minimum zero-sum uplift pricing approach that increases the price above the marginal cost and transfers all the additional payments from the profitable suppliers to the unprofitable suppliers to compensate their losses. Ruiz, Conejo, and Gabriel (2012) proposed a primal-dual approach to find the market prices that minimize the social welfare reduction due to the schedules inconsistencies and ensure non-negative profits for producers. Nevertheless, this approach does not guarantee the opportunity cost compensation for the generators, and the competitive equilibrium at the centralized dispatch solution is not achieved. Araoz and Jörnsten (2011) proposed a semi-Lagrangian relaxation approach to calculate a uniform market price that produces the same solution as the original centralized dispatch problem while ensuring no losses for the suppliers. A zero-sum uplift pricing scheme that minimizes the maximum contribution to the uplifts in the case of price-sensitive demand was developed by Van Vyve, (2011). Ring (1995), Gribik, Hogan, and Pope (2007), and Hogan and Ring (2003) developed the convex hull pricing method (CHP), which stays within the linear pricing framework with the uplift payments introduced to ensure the stability of the centralized dispatch solution. In this method, each market player is compensated the lost profit calculated as the difference between profit inferred from the market player decentralized dispatch solution and its profit in the centralized dispatch solution at a given market price. Thus, in this pricing approach, it is assumed that each market player has an opportunity to supply (consume) any power volume belonging to its private feasible set. The total lost profit (and hence, the total uplift needed to stabilize the centralized dispatch) equals the duality gap emerging after the Lagrangian relaxation procedure is applied to the power balance constraint (Ring, 1995; Gribik et. al. 2007; Hogan and Ring, 2003). CHP produces a market price (a set of prices) that minimizes the total uplift payment needed to compensate the market players for these foregone opportunities. Also, CHP minimizes the bound on the redistribution resulting from the near-optimal schedules (Eldridge et al. 2019). Solving the Lagrangian dual problem of the unit commitment problem to apply the CHP method can be computationally expensive. Hua and Baldick (2017) proposed a polynomially solvable primal formulation for the Lagrangian dual problem.

The uplift allocation may result in confiscation on the supply (demand) side and/or distortion of the market player bids. The latter may take place if uplift charges are allocated among producers (consumers) in a way that reduces (inflates) their revenues (expenses) but prevents confiscation. In this case, the producers have economic incentives to inflate the power output costs in their bids, while consumers are motivated to indicate reduced benefit from power consumption (provided that the output/consumption volumes cleared by the market are unchanged) to show close to zero profit obtained at the market price and avoid the uplift charge allocation. In addition, large uplift payments may result in market power abuse by the market players. Also, the uplifts decrease market pricing transparency and suppress economic signals. Therefore, it is critical to reduce the total uplift payment needed to support the centralized dispatch solution. An approach to achieve that based on the introduction of one redundant linear constraint in the form of an inequality and the associated price was proposed by Zhang et al. (2009). The resulting uplift is reduced compared to CHP at the expense of having a new service (a unit being online) and the corresponding price introduced at the market. The redundant nature of the new constraint provides that the feasible set of the primal problem is unchanged, while the linearity of the new constraint ensures that the duality gap is unaffected. In this approach, the total uplift payment equals the duality gap minus the nonnegative contribution from the redundant constraint.

Borokhov (2018) proposed identifying the opportunities available to the market players at the market and modifying CHP to compensate the lost profit only for the opportunities that are foregone by accepting the centralized dispatch solution. As a result, the total uplift payment needed to support the centralized dispatch solution is reduced (or unchanged) compared to that in CHP. However, the set of output (consumption) volumes corresponding to the opportunities available for a producer (consumer) is defined using the optimization problem and not constructed explicitly, which complicates the general analysis of the resulting pricing outcomes. In this paper, we apply the modified CHP to a single-period power market based on a uninode power system, i.e. a power system without transit losses, network and inter-temporal (such as ramping) constraints, with perfectly inelastic demand (fixed load) and zero minimum output limits of all generators and obtain explicit analytic expressions for the set of output volumes corresponding to the opportunities available for a producer at the market. This allows identifying a type of units that require special treatment in the proposed modified pricing approach compared to CHP and showing that the rest of the units may participate in the price setting procedure with their original private feasible set. We also establish that the market prices obtained from the modified CHP are "no lower" than the prices resulting from CHP.

The paper is organized as follows. We start with a short review of CHP in Section 2 and define a set of power outputs that can be supplied by generators to the centralized power market in Section 3. In Section 4, we formulate our proposal in terms of the modified CHP, show that in the convex case it is equivalent to standard marginal pricing, and identify the class of power systems for which CHP and the proposed method may result in the different sets of market prices. In Section 5, we analyze the structure of the set of market prices resulting from the modified CHP. Section 6 contains examples of power systems with comparisons of market prices and the associated total uplift payments obtained in CHP and the proposed method. The discussion and conclusion are presented in Section 7 and Section 8, respectively. The technical proofs are given in the Appendix.

## 2. The convex hull pricing

### 2.1 The centralized dispatch optimization problem

Consider a centrally dispatched single-period uninode power market with fixed demand $d$ and $|I|$ generating units bidding cost functions $C_i(x_i)$, $x_i = (u_i, g_i)$, $i \in I$, $I = \{1,...,|I|\}$, with the unit status binary variables $u_i$ taking values in the set $Z_2 = \{0,1\}$ (with 0 for an offline unit and 1 for an online unit) and output volumes $g_i$, $g_i \in R$, $0 \le g_i \le g_i^{\max}$, where $g_i^{\max} > 0$ denotes the unit maximum output limit. It is assumed that the system has sufficient generation capacity to meet demand: $d \le \sum_{i \in I} g_i^{\max}$. The generator cost functions are supposed to have a structure $C_i(x_i) = c_i(g_i) + w_i u_i$ with the fixed cost $w_i$, $w_i \ge 0$, and a non-decreasing convex continuous function $c_i(g_i)$ defined on $0 \le g_i \le g_i^{\max}$ with $c_i(0) = 0$. The fixed cost corresponds either to start-up and/or no-load cost. For simplicity, we will refer to it as the start-up cost. Initially, all the units are assumed to be offline. The centralized dispatch optimization problem (the primal problem) with the optimization (decision) variables $x = (x_1,..,x_n)$ is formulated as

$$\nu = \min_{\substack{x,\\ x_i \in X_i, \forall i \in I,\\ \sum_{i\in I} g_i = d}} \sum_{i\in I} C_i(x_i), \qquad (1)$$

where $\nu$ is the total cost to meet demand $d$ and $X_i$ denotes a set specified by the private constraints of a generator $i$: $X_i = \{x_i \mid u_i \in Z_2, g_i \in R, 0 \le g_i \le u_i g_i^{\max}\}$, $\forall i \in I$. (We note the following property of the given formulation of generator private constraints: zero output is possible for both online and offline statuses of the unit. If the start-up cost is nonzero, then the online status with zero output will not be an outcome of problem (1). If the start-up cost vanishes, then the outcome of (1) with unit zero output can have any of two possible unit statuses, which indicates the redundancy of the status variable as an optimization variable for the unit with zero start-up cost.) Let $\Omega$ denote the feasible set of the primal problem (1). Due to $d \le \sum_{i\in I} g_i^{\max}$, the set $\Omega$ is nonempty. Also, $\Omega$ is compact as the intersection of the compact sets $X_i$, $\forall i \in I$, and the closed set specified by the power balance constraint. Let $x^* = (x_1^*,..,x_n^*)$, $x_i^* = (u_i^*, g_i^*)$, $i \in I$, denote a solution to (1).

### 2.2 Decentralized dispatch optimization problem

Given a market price $p$, the decentralized dispatch problem for generating unit $i$ is formulated as

$$\pi_i(p) = \max_{x_i \in X_i} \pi_i(p, x_i), \text{ with } \pi_i(p, x_i) = p g_i - C_i(x_i), \quad (2)$$

and defines the supply curve of the unit $i$ as a set of points $(g_i, p)$ with $g_i \in \arg\max_{x_i \in X_i} \pi_i(p, x_i)$. To proceed further, we state some well-known mathematical facts about the optimization problem in question. We note that since $\pi_i(p)$ is a point-wise maximum of a function linear in $p$, it is convex in $p$ with the well-defined subdifferential $\partial\pi_i(p)$. Since $dom \quad \pi_i(p) = R$ is an open set, $\pi_i(p)$ is also continuous on $R$. Let us define the minimum economic output $g_i^{ec.\min}$ as follows: if $w_i = 0$, then $g_i^{ec.\min} = 0$; if $w_i \ne 0$, then $g_i^{ec.\min}$ is the lowest solution (if any) to the equation $[w_i + c_i(g_i^{ec.\min})]/g_i^{ec.\min} \in \partial c_i(g_i^{ec.\min})$ for $0 < g_i^{ec.\min} < g_i^{\max}$; if there is no solution, then $g_i^{ec.\min} = g_i^{\max}$. (We note that for $0 < g_i^{ec.\min} < g_i^{\max}$, the output $g_i^{ec.\min}$ can be equivalently defined as the lowest solution to $[w_i + c_i(g_i^{ec.\min})]/g_i^{ec.\min} \le \partial_+ c_i(g_i^{ec.\min})$, which is identical to the condition $\partial_+\{[w_i + c_i(g_i^{ec.\min})]/g_i^{ec.\min}\} \ge 0$, where $\partial_+$ denotes the right derivative.) Thus, $g_i^{ec.\min}$ depends only on the generator cost function and private feasible set $X_i$. For any given price $p$, the set of maximizers of (2) has the output volumes in the set $\{0\} \cup [g_i^{ec.\min}, g_i^{\max}]$. Thus, if the start-up cost is nonzero, then the supply curve has a gap as the curve has no points with output in the range $(0, g_i^{ec.\min})$. The output volumes from this range are never supplied in the decentralized dispatch problem under any $p$. (We note that $(0, g_i^{ec.\min}) \subset \partial\pi_i(p)$ for $p = [w_i + c_i(g_i^{ec.\min})]/g_i^{ec.\min}$.) For example, in the case of nonzero start-up cost and linear $c_i(g_i)$, we have $g_i^{ec.\min} = g_i^{\max}$ and the output volume that maximizes (2) equals either zero and/or $g_i^{\max}$ depending on the value of $p$, thus output volumes from the open interval $(0, g_i^{\max})$ do not maximize (2) at any market price $p$. However, for $p = [w_i + c_i(g_i^{\max})]/g_i^{\max}$ we have $\partial\pi_i(p) = [0, g_i^{\max}]$.

The optimization over the binary variable $u_i$ for a given value of $g_i$ allows to exclude the binary variable from (2) at the expense of having a discontinuity introduced in the cost function:

$$\pi_i(p) = \max_{g_i \in [0; g_i^{\max}]} \left[ p g_i - f_i(g_i) \right], \quad (3)$$

with $f_i(g_i) = w_i\theta(g_i) + c_i(g_i)$, $dom \quad f_i(g_i) = [0, g_i^{\max}]$, and the step-function $\theta(g_i)$ defined equal to 1 for $g_i > 0$ and 0 for $g_i \le 0$. Let us define $f_i(g_i) = +\infty$ outside $dom \quad f_i(g_i)$ and extend the feasible set in (3) to $R$, then $\pi_i(p)$ is the Fenchel convex conjugate of $f_i(g_i)$: $\pi_i = f_i^c$. An important property of the Fenchel conjugation (Rockafellar, 1970) is that $f_i\ (g_i)$ in (3) can be replaced by the greatest closed convex function majorized by $f_i\ (g_i)$, which we denote as $f_i^h(g_i)$. The function $f_i^h(g_i)$ is known as the closed convex hull of $f_i\ (g_i)$ and can be formally obtained by the double Fenchel conjugation of $f_i\ (g_i)$, (Rockafellar, 1970). We have $dom \quad f_i^h(g_i) = [0, g_i^{\max}]$. The function $f_i^h(g_i)$ is continuous on $(0, g_i^{\max})$ and takes infinite values on $(-\infty, 0) \cup (g_i^{\max}, +\infty)$. In the general case, this replacement results in a different set of maximizers for (3). Since (3) stays valid if $f_i(g_i)$ and the feasible set $[0, g_i^{\max}]$ are replaced by $f_i^h(g_i)$ and $R$ respectively, the application of the inversion rule for subgradient relations (Rockafellar and Wets, 1997) yields

$$g_i \in \partial\pi_i(p) \leftrightarrow p \in \partial f_i^h(g_i). \quad (4)$$

Thus, the set of output-price points of the form $(g_i, \partial f_i^h(g_i))$ can be equivalently represented as $(\partial\pi_i(p), p)$. This set of points would be the supply curve for the unit $i$ if its cost function were given by $f_i^h(g_i)$. Eq. (4) implies

$$\partial\pi_i(p) = \begin{cases} 0, p < [w_i + c_i(g_i^{ec.\min})]/g_i^{ec.\min} \\ [0, \gamma_{i\max}], p = [w_i + c_i(g_i^{ec.\min})]/g_i^{ec.\min} \\ \ldots, p > [w_i + c_i(g_i^{ec.\min})]/g_i^{ec.\min} \end{cases}, \quad (5)$$

where $\gamma_{i\max}$ is the maximum output volume $g_i$ satisfying $[w_i + c_i(g_i)]/g_i \in \partial c_i(g_i)$ for $0 < g_i \le g_i^{\max}$ and the dots denote the elements (which may depend on $p$) no lower than $\gamma_{i\max}$. We note that if $g_i^{ec.\min} = 0$, then (5) is well-defined because both $c_i(0) = 0$ and $g_i^{\max} > 0$ as well as the convexity of $c_i(g_i)$ imply that $\lim_{x_i \to +0} c_i(g_i)/g_i = \partial_+ c_i(0)$, which is finite.

**2.3 The dual problem**

The dual of the primal optimization problem (1) is formulated as:

$$v^D = \max_{p \in R} \min_{x_i \in X_i, \forall i \in I} [p(d - \sum\nolimits_{i \in I} g_i) + \sum\nolimits_{i \in I} C_i(x_i)] = \max_{p \in R} [pd - \sum\nolimits_{i \in I} \pi_i(p)]. \quad (6)$$

Since the RHS of (6) is an unconstrained maximization problem of a concave function, its set of maximizers, which we denote as $P^+$, is given by the solutions to $d \in \sum_{i \in I} \partial \pi_i(p)$. It is straightforward to see that $P^+$ is a nonempty closed convex set and contains nonnegative elements only. Thus, it has one of the following forms: a singleton $\{a\}$, a bounded closed interval of the form $[a, b]$, or a ray $[a, +\infty)$ with some $a, b \in R$, $0 \le a < b$. Let $x_i^+ = (u_i^+, g_i^+)$ denote a maximizer of the problem (2) with $p^+ \in P^+$. Define $\pi_i^*(p^+) = \pi_i(p^+, x_i^*)$, $\pi_i^+(p^+) = \pi_i(p^+, x_i^+)$. The duality gap is given by

$$v - v^D = \sum\nolimits_{i \in I} [\pi_i^+(p^+) - \pi_i^*(p^+)], \quad p^+ \in P^+, \quad (7)$$

and, according to CHP approach, represents the sum of the generator lost profits associated with opportunities to supply power output of the dispatch $x^+ = (x_1^+, .., x_n^+)$ foregone by accepting the dispatch $x^*$ at a price $p^+ \in P^+$. Thus, the duality gap equals the total uplift, and minimizing the total uplift is identical to solving the dual problem (6). Clearly, the duality gap is nonnegative. In the case under consideration, the non-zero duality gap may occur only if the market price in CHP method is unique, which is formalized by the following proposition.

*Proposition 1*: If $P^+$ is not a singleton, then $v = v^D$ and the duality gap is zero.

*Proof.* In the Appendix.

CHP implies the distribution of the amount (7) to the generators as the uplift payments to ensure that no generator acting as a price-taker (i.e. leaving aside issues related to the exercise of market power) has an economic incentive to change its output $g_i^*$ given the market price $p^+$. This can be viewed as utilizing the revenue function for the generator $i$ of the form

$$R_i(p^+, x_i) = p^+ g_i + \delta_{\{x_i^*\}}(x_i)[\pi_i^+(p^+) - \pi_i^*(p^+)], \quad (8)$$

with $p^+ \in P^+$ and the function $\delta_A(x_i)$, which is the 0-1 indicator function of a set $A$ : $\delta_A(x_i) = 1$ for $x_i \in A$ and $\delta_A(x_i) = 0$, otherwise. Thus, $\delta_{\{x_i^*\}}(x_i)$ is the 0-1 indicator function of a singleton $\{x_i^*\}$. (In practice, $\delta_{\{x_i^*\}}(x_i)$ can be defined equal to 1 if $| g_i - g_i^* |$ lies within a certain tolerance band.) An important attractive property of the prices resulting from CHP is that they are monotonically increasing in load since the subdifferential of a convex function is a non-decreasing operator; however, due to the uplift payments, the aggregate generator revenue is generally not a monotonous function of load (Gribik, Hogan, and Pope, 2007). A price $p$ is said to support the solution $x^*$ if $x_i^* \in \arg\max_{x_i \in X_i} [pg_i - C_i(x_i)]$, $\forall i \in I$. The price that supports a solution to the primal problem exists if and only if the duality gap is zero. In this case, the set of prices supporting the solution $x^*$ is independent of the choice of $x^*$ (if the primal problem has multiple solutions) and is given by the set of maximizers for (6). The aggregate supply at a given price is the sum of generator outputs obtained by solving the corresponding decentralized dispatch problems at this price. The dual problem provides a framework to find a market price (a set of prices) that corresponds to a transition from a shortage of aggregate supply compared to demand $d$ to a surplus of aggregate supply over $d$. If the aggregate supply and demand curves intersect at some $p$, then the power balance constraint holds at the intersection segment, and $x^+$ (some $x^+$ if (2) has multiple solutions) belongs to the feasible set of the primal problem. Therefore, the duality gap is zero and $p \in P^+$. (Here we used the fact that the power balance constraint relaxed in the dual problem formulation has the form of equality. If an inequality constraint is relaxed and some optimal point of the relaxed problem is feasible in the primal problem, then the duality gap could still be present since the complementarity condition may not hold.) The converse is also true: if the duality gap is zero, then any $p \in P^+$ supports $x^*$ and the aggregate supply and demand curves intersect at the price $p$ because $x^*$ satisfies the power balance constraint. Therefore, a duality gap in the model under consideration occurs only if the aggregate supply and demand curves do not intersect. Thus, in the case of a uninode single-period power system under consideration, CHP produces a market price that corresponds to either an intersection of aggregate supply and demand curves (zero duality gap) or a transition from undersupply, i.e. demand exceeds aggregate supply, to oversupply, i.e. aggregate supply exceeds demand, (nonzero duality gap), Fig. 1. Note that (3) provides a straightforward way to find $P^+$ for a uninode single-period power system since it is the set of the marginal prices in the convex centralized dispatch problem obtained from the primal problem when the generator cost functions $C_i(x_i)$ are replaced by $f_i^h(g_i)$. Graphically, $P^+$ corresponds to the intersection of the demand curve with the new aggregate supply curve constructed from the convex cost functions $f_i^h(g_i)$, $i \in \boldsymbol{I}$. To illustrate some implications of CHP, we consider the following example. Schiro et al. (2016) made a comprehensive study of important properties of CHP.

*Example 1.* Consider a power system with fixed demand $d$ and a single supplier having the cost function $C(x) = wu + ag$ with start-up cost $W$, $w > 0$, and the constant marginal cost $a$, $a \geq 0$, zero minimum output limit, and the maximum output limit $g^{\max}$, which exceeds demand, $d < g^{\max}$. Clearly, the primal problem solution yields $u^* = 1$, $g^* = d$, and the marginal price $a$. The supply curve reconstructed from the generator decentralized dispatch problem is comprised of two disjoint segments $g = 0$ for $p \leq a + w / g^{\max}$ and $g = g^{\max}$ for $p \geq a + w / g^{\max}$ and is, therefore, discontinuous. We also observe that the supply curve has no point with the output equal demand, and, as a result, the supply and demand curves do not intersect. Furthermore, since $g^{ec.\min} = g^{\max}$, the generator output $g^*$ is below $g^{ec.\min}$. The application of CHP yields the market price $p^+ = a + w / g^{\max}$ (the set $P^+$ is singleton), which is below the generator average cost for the output $d$. The price $p^+$ implies two possibilities for the unit status-output variables in the decentralized dispatch problem: $u^+ = g^+ = 0$, and $u^+ = 1$, $g^+ = g^{\max}$, both yielding $\pi(p^+, x^+) = 0$. Hence, the generator is compensated with an uplift payment of $w(1 - d / g^{\max})$, which results in zero generator profit for output $d$. We note that if the price were set above $p^+$, then according to the CHP principle, the generator would have to be compensated for the foregone opportunity to profitably supply $g^{\max}$, which increases the uplift payment. However, the generator is not able to supply any nonzero volume other than $d$. Therefore, the generator had not lost any opportunity to supply any higher output by accepting the centralized dispatch solution. Also, since the price $p^+$ is below the generator average output cost, it may deter new potential entry (able to fully replace the incumbent generator) from entering the market because $p^+$ alone, without the uplift payment, underestimates the level of the average output cost needed to replace the producer. Thus, $g^{\max}$ -dependence of the market price seems counterintuitive. It is more desirable to have the price that is independent of the infeasible output volumes, i.e. volumes above $d$. This is achieved, for example, when the price is set to $(a + w/d)$, in which case no uplift payment is needed at all if the infeasible output volumes are removed from the lost profit calculation. The reason for $g^{\max}$ -dependence of $p^+$ is that $g^{ec.\min} = g^{\max}$ with $g^{\max} > d$. As a result, one possible way to have the dual problem optimal output closer to $d$ in case of $u^+ = 1$ is to replace $g^{\max}$ by $d + \varepsilon$ with some small $\varepsilon > 0$ in the dual problem. These observations lead to a principal question of whether the output volumes in the range $(d, g^{max}]$ should be used in the lost profit calculation. We address this issue in Section 4.

The compensation of the profit lost due to the foregone opportunity implies that the market player could receive this additional profit if it were not for the centralized dispatch. The example above illustrates that some of the economic opportunities, which are treated as foregone and compensated in CHP, may not be realizable by the market players in the absence of the centrally coordinated market. Hence, the output opportunities available to a generator in the power system should be examined in more detail.

## 3. The power outputs attainable at the market

The lost profit compensations associated with foregone opportunities are needed to ensure the economic stability of the centralized dispatch outcome. This implies that in the absence of the centralized dispatch each generator has a set of possible output volumes that it can supply to the power system provided that these volumes can be accepted from both technological and economic considerations (we refer to this set as a set of technologically and economically feasible output volumes for a given producer). We determine such a set for each generator as follows. Let us allow a generator $i_0$ to engage in the bilateral contracts for power with the other market participants (both generators and consumers) by paying to the other generators the full cost of contracted output volumes and receiving payments from consumers for the contracted volumes in the amount of their benefit from power consumption. The costs/benefits in the described procedure are calculated using the bids submitted by the producers/consumers to the centrally coordinated market. These bilateral contracts are the financial contracts, which are accounted for at the financial settlement stage only, and the market participants with bilateral contracts still submit their bids to the centralized market. Since demand is fixed, we formally require all demand volumes to be fully contracted, which leaves the feasible set of the optimization problem unchanged. We also require the resulting output/consumption schedule to be attainable as a centralized dispatch outcome. In this setting, only the market participants with the contracts are allowed to participate in the centralized dispatch optimization problem using their original bids with volumes restricted by the contracted output/consumption volumes (the generator $i_0$ also submits the bid for its contracted output, which is equal to the netted contracted volume). To ensure economic feasibility of the resulting schedule, we require the total financial effect for the generator $i_0$ from all the contracts to be nonnegative. This can be formalized as follows. We say that a generator $i_0$ has an opportunity to supply output volume $\bar{g}_{i_0}$, $0 \le \bar{g}_{i_0} \le g_{i_0}^{max}$, at the centralized dispatch market if there exists a set of $\bar{g}_i$, $0 \le \bar{g}_i \le g_i^{max}$, $i \in I \setminus \{i_0\}$, and a set of the corresponding statuses $\bar{u}_i$, $i \in I$, so that $\bar{x} = (\bar{x}_1, .., \bar{x}_n)$ with $\bar{x}_i = (\bar{u}_i, \bar{g}_i)$ satisfies the following condition:

$$\bar{x} \in \arg \min_{\substack{x, \\ x_i \in X_i, \forall i \in I, \\ x \le \bar{x}, \\ \sum_{i \in I} g_i = d}} \sum_{i \in I} C_i(x_i) . \qquad (9)$$

Thus, the feasible set in (9) is given by that of (1) supplemented with $u_i \le \bar{u}_i$ , $g_i \le \bar{g}_i$ , $\forall i \in I$ . Clearly, we have $\bar{x} \in \Omega$. Let us denote by $\bar{\Omega}$ the set of all $\bar{x}$ satisfying (9). Since the generator $i_0$, aside from producing output $\bar{g}_{i_0}$ , merely acts as an intermediary (retaining the market surplus), the set $\bar{\Omega}$ is independent of the choice of $i_0$. To make a transition from $\bar{\Omega}$ to a set of status-output volumes for a given producer $i'$, let us denote as $\bar{\Omega}_{i'}$ a set containing all $x_{i'}$ such that there exist a collection of $x_i$ , $i \in I$, $i \ne i'$, so that $x \in \bar{\Omega}$ . The set $\bar{\Omega}_{i'}$ can be viewed as projection of $\bar{\Omega}$ on the set $X_{i'}$. We have $\bar{\Omega}_i \subset X_i$, $\bar{\Omega} \subset \times_{i \in I} \bar{\Omega}_i$ , but in the general case $\bar{\Omega} \ne \times_{i \in I} \bar{\Omega}_i$ . Also, the primal problem outcome $x^*$ can be realized through a set of bilateral contracts since $x = x^*$ satisfies (9). Therefore, $x^* \in \bar{\Omega}$, $x_i^* \in \bar{\Omega}_i$ .

Let us denote as $N$ all elements $x$ from $\Omega$ with $u_i = 1$, $g_i = 0$ for at least one $i$ such that $w_i \ne 0$ (if there are no such elements, then $N$ is an empty set). Now we show that all elements of $\Omega' = \Omega \setminus N$ can be realized using some set of the bilateral contracts.

*Proposition 2*: $\Omega' = \bar{\Omega}$.

*Proof*. Clearly, we have $\bar{\Omega} \subset \Omega'$, so we need to show that $\Omega' \subset \bar{\Omega}$. Let $x' \subset \Omega'$. Since the corresponding $g'$ satisfies $\sum_{i \in I} g_i' = d$ and $0 \le g_i' \le g_i^{\max}$ , $\forall i \in I$, we conclude that $x'$ satisfies (9) for $\bar{x} = x'$ because $x = x'$ is the only feasible point in the optimization problem (9). Thus, we have $\Omega' \subset \bar{\Omega}$, which implies $\Omega' = \bar{\Omega}$. ■

Proposition 2 allows us to find $\bar{\Omega}_{i_0}$ explicitly by projecting $\Omega'$ on $X_{i_0}$: the corresponding values of the output volumes $\bar{g}_{i_0}$ are given by the closed interval $[\bar{g}_{i_0}^{\min}, \bar{g}_{i_0}^{\max}]$ with $\bar{g}_{i_0}^{\min} = \max \;\; (d - \sum_{i: i \in I, i \ne i_0} g_i^{\max}; 0)$, $\bar{g}_{i_0}^{\max} = \min \;\; (d; g_{i_0}^{\max})$. (We note that this expression for possible values of $\bar{g}_{i_0}$ heavily relies on the assumptions of fixed load and zero minimum output limits of all generators in the power system.) Thus, if for a unit $i$ both zero generator output volume and maximum generator output volume are feasible in the primal problem, then for this unit the set of technologically and economically feasible output volumes coincides with that specified by its private feasible set, i.e. $[0, g_i^{\max}]$. If $N$ is nonempty, then the sets $\Omega'$ and $\bar{\Omega}$ are noncompact. Hence, the set $\overline{\Omega}_{i_0}$ is noncompact iff both $w_{i_0} \neq 0$ and $\overline{g}_{i_0}^{\min} = 0$. In this case, we have $\overline{\Omega}_{i_0} = \{(0,0)\} \cup \{(u_{i_0}, g_{i_0}) \mid u_{i_0} = 1, 0 < g_{i_0} \leq \overline{g}_{i_0}^{\max}\}$.

## 4. The modified convex hull pricing

### 4.1 The modified primal problem and the corresponding dual problem

Let us define a modified primal problem

$$\tilde{\nu}(\varepsilon) = \min_{\substack{x, \\ x_i \in \tilde{X}_i(\varepsilon), \forall i \in I \\ \sum_{i \in I} g_i = d}} \sum_{i \in I} C_i(x_i) \qquad (10)$$

with the modified private feasible set $\tilde{X}_i(\varepsilon_i)$ of a generator $i$ given by

$$\tilde{X}_i(\varepsilon_i) = \{(0,0)\} \cup \bigcup_{\overline{x}_i \in \overline{\Omega}_i} \Delta_i(\overline{x}_i, \varepsilon_i), \text{ where } \Delta_i(\overline{x}_i, \varepsilon_i) = \{(u_i, g_i) \mid (u_i, g_i) \in X_i, \left| g_i - \overline{g}_i \right| \leq \varepsilon_i\}, \quad (11)$$

with some $\varepsilon_i > 0$, $\varepsilon = (\varepsilon_1, .., \varepsilon_n)$. The inclusion of $\Delta_i(\overline{x}_i, \varepsilon_i)$ for each element $\overline{x}_i \in \overline{\Omega}_i$ is needed to indicate in the dual problem whether at a given price a generator $i$ is willing to supply some more/less power than $\overline{g}_i$ within its private feasible set $X_i$. Thus, for each $\overline{g}_i$ with $\overline{x}_i \in \overline{\Omega}_i$, the sets $\tilde{X}_i(\varepsilon_i)$ and $X_i$ have identical output volumes in some closed neighborhood of $\overline{g}_i$. Regarding the need to ensure that $\tilde{X}_i(\varepsilon_i)$ includes a point $(u_i = 0, g_i = 0)$, which implies zero generator profit at any market price, we have the following comment. If $\overline{\Omega}_i$ includes only the elements with the unit $i$ having output no lower than some positive value, then for sufficiently small $\varepsilon_i$ all elements of $\bigcup_{\overline{x}_i \in \overline{\Omega}_i} \Delta_i(\overline{x}_i, \varepsilon_i)$ correspond to the nonzero output of the unit. In this case, if $(u_i = 0, g_i = 0)$ is not included in $\tilde{X}_i(\varepsilon_i)$, then for a given price in the decentralized dispatch problem the unit will find its optimal output volume disregarding the start-up cost $w_i$. Therefore, in such a case, $w_i$ will contribute neither to the market price nor the uplift payment for the unit, which may result in the confiscatory pricing for the generator. (Note that if $w_i = 0$, then zero generator profit at any market price is also realized at $(u_i = 1, g_i = 0)$, and any of these elements can be included in $\tilde{X}_i(\varepsilon_i)$ to guarantee the non-confiscatory pricing for power.) By construction, we have $\tilde{X}_i(\varepsilon_i) \subset X_i$. Let us denote by $\tilde{\Omega}(\varepsilon)$ the feasible set of the modified primal problem (10). Since $x^* \in \overline{\Omega} \subset \tilde{\Omega}(\varepsilon)$ for any $x^*$ and $\tilde{\Omega}(\varepsilon) \subset \Omega$, we conclude that $\nu = \tilde{\nu}(\varepsilon)$ and both the primal and the modified primal problems have identical sets of the maximizers. The analysis of Section 3 yields the following expressions for $\tilde{X}_i(\varepsilon_i)$, $\forall \varepsilon_i > 0$:

- if $\varepsilon_i < \overline{g}_i^{\min}$, then $\tilde{X}_i(\varepsilon_i) = \{(0,0)\} \cup \{(u_i, g_i) \mid u_i = 1, \overline{g}_i^{\min} - \varepsilon_i \le g_i \le \min[\ \overline{g}_i^{\max} + \varepsilon_i; g_i^{\max}\ ]\}$;
- if $\varepsilon_i \ge \overline{g}_i^{\min}$, then $\tilde{X}_i(\varepsilon_i) = \{(0,0)\} \cup \{(u_i, g_i) \mid u_i = 1, 0 \le g_i \le \min[\ \overline{g}_i^{\max} + \varepsilon_i; g_i^{\max}\ ]\}$;

which illustrate a need to introduce $\varepsilon_i$. For a given market price, consider the decentralized dispatch problem for a unit $i$ with the feasible set $\tilde{X}_i(\varepsilon_i)$. If $\varepsilon_i = 0$ and $\overline{g}_i^{\min} > 0$ is the optimal unit $i$ output, then the generator either sells all output volumes at a price no lower than its marginal cost of output or the generator makes nonnegative profit but sells some part of the output volumes below its marginal cost. In the latter case, it has economic incentives to decrease its output. Likewise, if $\overline{g}_i^{\max}$ is the optimal unit output and $\varepsilon_i = 0$, then either $\overline{g}_i^{\max} = g_i^{\max}$ and the unit will not change its output if the market price increases or the unit will increase its output if the price increases. The introduction of $\varepsilon_i > 0$ allows differentiating between these cases. We also note that $\tilde{X}_i(\varepsilon_i)$ is compact, $\forall \varepsilon_i > 0$, $\forall i \in I$. The dual of the modified primal problem has the form

$$\tilde{v}^D(\varepsilon) = \max_{p \in R} \left( pd - \sum_{i \in I} \tilde{\pi}_i(p, \varepsilon_i) \right), \text{ with } \tilde{\pi}_i(p, \varepsilon_i) = \max_{x_i \in \tilde{X}_i(\varepsilon_i)} \pi_i(p, x_i). \quad (12)$$

According to (11), the set $\tilde{X}_i(\varepsilon_i)$ for a unit with $w_i \ne 0$ may include the point $(u_i = 1, g_i = 0)$. If (11) is modified to exclude that point from $\Delta_i(\overline{x}_i, \varepsilon_i)$ for such a unit, the resulting set $\tilde{X}_i(\varepsilon_i)$ becomes noncompact. However, since for a unit with $w_i \ne 0$ the point $(u_i = 1, g_i = 0)$ neither belongs to the minimizer of the (modified) primal problem nor maximizes $\tilde{\pi}_i(p, \varepsilon_i)$, the inclusion of that point in $\Delta_i(\overline{x}_i, \varepsilon_i)$ affects none of $\tilde{v}(\varepsilon)$, $\tilde{v}^D(\varepsilon)$, and $\tilde{\pi}_i(p, \varepsilon_i)$. Thus, we conclude that such a modification of $\Delta_i(\overline{x}_i, \varepsilon_i)$ changes neither the set of market prices obtained from (12) nor the individual generator uplifts. Likewise, the pricing outcomes are not affected if the point $(u_i = 0, g_i = 0)$ in the definition of $\tilde{X}_i(\varepsilon_i)$ is substituted or supplemented by $(u_i = 1, g_i = 0)$ for a generator with $w_i = 0$.

The function $\tilde{\pi}_i(p, \varepsilon_i)$ is convex in $p$ since it is a point-wise maximum of the function linear in $p$. We denote its subdifferential with respect to $p$ as $\partial \tilde{\pi}_i(p, \varepsilon_i)$. Let $\tilde{P}^+(\varepsilon)$ be a set of maximizers of (12), i.e. a set of prices that solve $d \in \sum_{i \in I} \partial \tilde{\pi}_i(p, \varepsilon_i)$. It is straightforward to verify that $\tilde{P}^+(\varepsilon)$ is nonempty. The relation $\tilde{X}_i(\varepsilon_i) \subset X_i$ implies $v^D \le \tilde{v}^D(\varepsilon)$. Hence, for $\forall \varepsilon_i > 0$, $\forall i \in I$, we have $v^D \le \tilde{v}^D(\varepsilon) \le \tilde{v}(\varepsilon) = v$, which entails that the duality gap in the modified problem is not higher than the duality gap in the original problem:

$$0 \le \tilde{v}(\varepsilon) - \tilde{v}^D(\varepsilon) \le v - v^D . \qquad (13)$$

Hence, if the duality gap in the original problem is zero, then the duality gap in the modified problem is also zero. (The converse is generally not true as illustrated in Example 3 below.) It can be shown that if $\tilde{P}^+(\varepsilon)$ is not a singleton, then $\tilde{v}^D(\varepsilon) = \tilde{v}(\varepsilon) = v$ and the duality gap vanishes (the proof is fully analogous to that of Proposition 1).

### 4.2 The market pricing proposal

The proposed modified CHP is formulated as follows. We suggest calculating the set $\tilde{P}^+(\varepsilon)$ in the limit as $\varepsilon \to +0$ and utilizing the resulting $\tilde{P}^+(+0)$ as the set of market prices. Let us define $\tilde{\pi}_i^+(p) = \tilde{\pi}_i(p,+0)$, $\forall i \in I$. To ensure the economic stability of the centralized dispatch outcome at the price $\tilde{p}^+ \in \tilde{P}^+(+0)$, we propose paying the generator $i$ the uplift payment of $[\tilde{\pi}_i^+(\tilde{p}^+) - \tilde{\pi}_i^*(\tilde{p}^+)]$ if it follows the centralized dispatch and penalizing the generator in the amount of $max[\pi_i \ (\tilde{p}^+, x_i) - \tilde{\pi}_i^+(\tilde{p}^+); 0]$ if its status-output is outside $\tilde{X}_i(+0)$. From $\tilde{v}(\varepsilon) - \tilde{v}^D(\varepsilon) = \sum_{i \in I} [\tilde{\pi}_i(p, \varepsilon_i) - \pi_i^*(p)]$ for $p \in \tilde{P}^+(\varepsilon)$, (7), and (13), we conclude that the total uplift needed to support the centralized dispatch solution at any price $\tilde{p}^+ \in \tilde{P}^+(+0)$ in the modified optimization problem (10) is no higher than that for the original problem (1) at any price $p^+ \in P^+$. Clearly, if $\tilde{P}^+(+0)$ contains more than one element, the individual generator uplifts are independent of the choice of $\tilde{p}^+ \in \tilde{P}^+(+0)$.

The proposed approach can be viewed as using the following expression for the generator revenue function for a given market price $p \in R$:

$$\tilde{R}_i(p, x_i) = pg_i + \delta_{\{x_i^*\}}(x_i)[\tilde{\pi}_i^+(p) - \pi_i^*(p)] - \delta_{X_i \setminus \tilde{X}_i(+0)}(x_i) \max[\pi_i \ (p, x_i) - \tilde{\pi}_i^+(p); 0], \ \forall x_i \in X_i, \forall i \in I .$$

Since $\max[\pi_i \ (p, x_i) - \tilde{\pi}_i^+(p); 0] = 0$, $\forall x_i \in \tilde{X}_i(+0)$, $\forall i \in I$, the revenue function $\tilde{R}_i(p, x_i)$ can be equivalently expressed as

$$\tilde{R}_i(p, x_i) = pg_i + \delta_{\{x_i^*\}}(x_i)[\tilde{\pi}_i^+(p) - \pi_i^*(p)] - \max[\pi_i \ (p, x_i) - \tilde{\pi}_i^+(p); 0], \ \forall x_i \in X_i, \forall i \in I . \quad (14)$$

For $\forall p \in R$, we also have

$$\begin{aligned} &\max_{x_i \in \tilde{X}_i(\varepsilon_i)} \pi_i(p, x_i) = \max_{x_i \in \tilde{X}_i(\varepsilon_i)} [pg_i - C_i(x_i)] = \max_{x_i \in \tilde{X}_i(\varepsilon_i)} [\tilde{R}_i(p, x_i) - C_i(x_i)] = \\ &\max\{ \max_{x_i \in \tilde{X}_i(\varepsilon_i)} [\tilde{R}_i(p, x_i) - C_i(x_i)]; \max_{x_i \in X_i \setminus \tilde{X}_i(\varepsilon_i)} [\tilde{R}_i(p, x_i) - C_i(x_i)\} = \max_{x_i \in X_i} [\tilde{R}_i(p, x_i) - C_i(x_i)] \end{aligned} .$$

Thus, given the market price $p \in R$, the decentralized dispatch problem of a generator in the proposed modified CHP is identical to the decentralized dispatch problem of the generator with the feasible set $X_i$ and the generator revenue function $\tilde{R}_i(p, x_i)$.

For the case of zero generator start-up costs ($w_i = 0$, $\forall i \in I$), both the duality gaps in the original and the modified problems are zero as these optimization problems become convex after the exclusion of the binary status variables and the total uplift vanishes. In this case, the set of prices obtained in CHP is identical to the set of prices given by the modified CHP, i.e. the marginal prices, (Borokhov, 2018). Moreover, a set of market prices resulting from the modified convex pricing approach is independent of $\varepsilon_i$. Thus, if the primal problem (1) is convex, then the infeasible output volumes of the unit $i$ (i.e. $[0, \bar{g}_i^{\min} - \varepsilon_i)$ if $\bar{g}_i^{\min} > 0$, and $(\bar{g}_i^{\max} + \varepsilon_i, g_i^{\max}]$ if $\bar{g}_i^{\max} < g_i^{\max}$, for sufficiently small $\varepsilon_i > 0$) can be removed from the generator private feasible set and, consequently, from the lost profit calculation with no effect on the set of market prices.

### 4.3 Utilization of the original private feasible sets in the proposed pricing method

For the non-convex case, we observe that the output $\bar{g}_i^{\min}$, $\forall i \in I$, can be formally set to zero in the dual of the modified primal problem: extending $\tilde{X}_i(\varepsilon_i)$ to include all the elements of $X_i$ with the output volumes in the range $[0, \bar{g}_i^{\max} + \varepsilon_i]$ affects neither the set of market prices nor the individual uplift payment received by the generators in the modified CHP. Define

$$\hat{v}^D(\varepsilon) = \max_{p \in R} \quad \left(pd - \sum_{i \in I} \hat{\pi}_i(p, \varepsilon_i)\right) \quad (15)$$

with $\hat{\pi}_i(p, \varepsilon_i) = \max_{x_i \in \hat{X}_i(\varepsilon_i)} \pi_i(p, X_i)$, $\hat{X}_i(\varepsilon_i) = \{(u_i, g_i) \mid (u_i, g_i) \in X_i, g_i \le \bar{g}_i^{\max} + \varepsilon_i\}$, $\forall i \in I$. Let $\hat{P}^+(\varepsilon)$ denote the set of maximizers of (15).

*Proposition* 3: For the optimization problems (12) and (15) with $\varepsilon_i > 0$, $\forall i \in I$, we have

- $\hat{P}^+(\varepsilon) = \tilde{P}^+(\varepsilon)$;
- $\hat{\pi}_i(p, \varepsilon_i) = \tilde{\pi}_i(p, \varepsilon_i)$, $\forall p \in \tilde{P}^+(\varepsilon), \forall i \in I$;
- $\hat{v}^D(\varepsilon) = \tilde{v}^D(\varepsilon)$.

*Proof.* In the Appendix.

Another view on Proposition 3 is the following. If $\bar{g}_i^{\min} \le \varepsilon_i$, then $\tilde{X}_i(\varepsilon_i) = \hat{X}_i(\varepsilon_i)$, $\tilde{\pi}_i(p, \varepsilon_i) = \hat{\pi}_i(p, \varepsilon_i)$, $\forall p \in R$. If $0 < \bar{g}_i^{\min} - \varepsilon_i \le \min(\bar{g}_i^{\max} + \varepsilon_i; g_i^{ec.\min})$, then since the output volumes from the open interval $(0, \min(\overline{g}_i^{\max} + \varepsilon_i; g_i^{ec.\min}))$ never maximize either $\tilde{\pi}_i(p,\varepsilon_i)$ or $\hat{\pi}_i(p,\varepsilon_i)$, we conclude that $\tilde{\pi}_i(p,\varepsilon_i) = \hat{\pi}_i(p,\varepsilon_i)$, $\forall p \in R$. At last, if $\overline{g}_i^{\min} - \varepsilon_i > g_i^{ec.\min}$, then Proposition 3 implies that the values of $p$ with $\tilde{\pi}_i(p,\varepsilon_i) \neq \hat{\pi}_i(p,\varepsilon_i)$ do not belong to $\hat{P}^+(\varepsilon) = \tilde{P}^+(\varepsilon)$.

When calculating the output volumes attainable at the centralized dispatch for a generator, we imposed the requirement that perfectly inelastic demand should be fully contracted, which implies that the power balance constraint of the primal problem holds. An alternative approach is to replace the fixed demand by the consumer bids with some benefit functions, e.g. with constant marginal benefit, and consider the limit as the marginal benefit goes to infinity. In this setting, some elements of the resulting set may not belong to the feasible set of the primal problem. For the sufficiently high marginal benefits of the consumers, this approach results in a set of technologically and economically feasible output volumes of the form $[0, \overline{g}_{i_0}^{\max}]$ for a unit $i_0$ with $\overline{g}_{i_0}^{\max}$ obtained in Section 3. Thus, in this case, $\overline{g}_{i_0}^{\min} = 0$ for any generator $i_0$. The statement of Proposition 3 entails that this alternative approach produces the same set of the market prices, generator profits, and the uplift payments.

Having shown that $\forall i \in I$ the output $\overline{g}_i^{\min}$ can be excluded from the consideration (i.e. formally set to zero) in the modified primal problem and its dual, we turn to the cases when $\overline{g}_i^{\max}$ can be disregarded as well. Let us partition $I$ into $\breve{I} = \{i \mid i \in I, g_i^{ec.\min} \leq \overline{g}_i^{\max}\}$ and $\hat{I} = \{i \mid i \in I, g_i^{ec.\min} > \overline{g}_i^{\max}\}$. Clearly, if $g_i^{\max} \leq d$, then $i \in \breve{I}$. In the case of $g_i^{\max} > d$ we have $i \in \breve{I}$ if $g_i^{ec.\min} \leq d$ and $i \in \hat{I}$ if $g_i^{ec.\min} > d$. Consider

$$\breve{v}^D(\varepsilon_{\hat{I}}) = \max_{p \in R} \left( pd - \sum_{i \in \breve{I}} \pi_i(p) - \sum_{i \in \hat{I}} \hat{\pi}_i(p,\varepsilon_i) \right), \quad (16)$$

with $\varepsilon_{\hat{I}} = \{\varepsilon_{\hat{i}}\}$, $\hat{i} \in \hat{I}$. Let $\breve{P}^+(\varepsilon_{\hat{I}})$ denote a set of maximizers of (16).

*Proposition* 4: For optimization problems (15) and (16) for $\varepsilon_i > 0$, $\forall i \in I$, we have

- $\breve{P}^+(\varepsilon_{\hat{I}}) = \hat{P}^+(\varepsilon)$;
- $\pi_i(p) = \hat{\pi}_i(p,\varepsilon_i)$, $\forall p \in \hat{P}^+(\varepsilon)$, $\forall i \in \breve{I}$;
- $\breve{v}^D(\varepsilon_{\hat{I}}) = \hat{v}^D(\varepsilon)$.

*Proof.* In the Appendix.

Thus, Propositions 3 and 4 imply that for $\varepsilon_i > 0$, $\forall i \in I$, we have $\tilde{P}^+(\varepsilon) = \hat{P}^+(\varepsilon) = \breve{P}^+(\varepsilon_{\hat{I}})$, which are independent of $\varepsilon_i$, $\forall i \in \breve{I}$, (but may depend on $\varepsilon_i$, $i \in \hat{I}$). Also, $\forall p \in \tilde{P}^+(+0)$, the expression (14) for the revenue function of any unit from the set $\breve{I}$ is equivalent to (8). Thus, (14) is needed only for $i \in \hat{I}$. Since $g_i^{ec.\min} > \overline{g}_i^{\max}$ implies $g_i^{ec.\min} > d$ (which, in turn, entails $g_i^{\max} > d$), each unit from the set $\hat{I}$ has the physical capability to meet the demand $d$. Thus, any unit from the set $\hat{I}$ exhibits natural monopoly behavior (has decreasing average cost function) for the output volumes in $(0, d]$. In what follows, we will refer to such a unit as the large natural monopoly generating unit (LNMGU). Thus, if no LNMGU is present in the power system, then both CHP and the modified CHP result in the identical sets of market prices and generator uplift payments. Consequently, these pricing methods may produce different outcomes only if the system has at least one LNMGU. (We emphasize that this conclusion for the considered uninode single-period power market with fixed load is tied to the assumption of vanishing minimum output limits of all generating units. If this assumption is relaxed, the modified CHP may result in the different pricing outcome compared to CHP even in the absence of LNMGU (Borokhov, 2018).) Clearly, each LNMGU has a nonzero start-up cost. Throughout the rest of the paper, we choose sufficiently small $\varepsilon_i$ so that $0 < \varepsilon_i < g_i^{ec.\min} - d$, $\forall i \in \hat{I}$. This ensures a decreasing average cost function for LNMGU output in the range $(0, d + \varepsilon_i]$. It can be shown that no more than one LNMGU can be online according to the primal problem solution (this is also true for the modified primal problem (10) as well). Now, we study LNMGUs in the context of the dual problem. Let us define

$$\hat{I}_{\min}(\varepsilon_{\hat{I}}) = \arg \min_{i \in \hat{I}} \; [w_i + c_i(d + \varepsilon_i)]/(d + \varepsilon_i).$$

*Proposition 5*: If $\hat{I} \neq \emptyset$, then for sufficiently small positive values of $\{\varepsilon_i\}$, $i \in \hat{I}$, removal of any group of LNMGUs from the consideration in the dual problem (16), provided that at least one LNMGU from $\hat{I}_{\min}(\varepsilon_{\hat{I}})$ remains, will not change the resulting set of market prices.

*Proof*. In the Appendix.

Thus, all LNMGUs except for any one belonging to $\hat{I}_{\min}(\varepsilon_{\hat{I}})$ for the given values of $\varepsilon_i > 0$, $i \in \hat{I}$, can be removed from the consideration in the dual problem (16) without affecting both the set of market prices $\tilde{P}^{+}(\varepsilon)$ and the individual uplifts of all generating units. We note that the statement of Proposition 5 is also valid for the dual problems (6), (12), (15) with possibly different relevant LNMGUs. Note that the set $\hat{I}_{\min}(\varepsilon_{\hat{I}})$ may depend on values of $\{\varepsilon_i\}$, $i \in \hat{I}$. Thus, at most two LNMGUs are relevant: one for the primal problem solution for the power output and the other for the dual problem solution for the market prices. Also, the online LNMGU in the solution for the primal problem (1), which is also a solution for the modified primal problem (10), and the LNMGU setting the price in (16) might be different even if a solution to the primal problem is unique. This possibility is illustrated in the example below.

*Example 2.* Consider a single-period uninode power system with fixed demand $d$ and three generating units with the cost functions given by: $C_1(x_1)=0$, $0 \le g_1 \le d/4$; $C_2(x_2)=d^2u_2+g_2^2/2$, $0 \le g_2 \le 2d$ ; $C_3(x_3)=1.4d^2u_3$, $0 \le g_3 \le 2d$. Thus, units 1 and 3 have zero variable cost of output. It is straightforward to verify that $g_1^{ec.\min}=0$, $g_2^{ec.\min}=\sqrt{2}d$, $g_3^{ec.\min}=2d$. Therefore, units 2 and 3 are LNMGUs. Both the primal problem (1) and the modified primal problem (10) have a unique optimal point with unit 1 having output $d/4$, unit 2 producing $3d/4$, and unit 3 being offline. Thus, unit 2 is the online LNMGU in the unique primal problem solution. The solution for the problem (16) has unit 1 producing $d/4$, while the rest of the demand is satisfied with either of unit 2 or unit 3 depending on which unit has the lowest average total cost for the output volumes $d+\varepsilon_2$ and $d+\varepsilon_3$, respectively: $(d^2+(d+\varepsilon_2)^2/2)/(d+\varepsilon_2)$ for unit 2 and $1.4d^2/(d+\varepsilon_3)$ for unit 3. For sufficiently low $\varepsilon_2$, we have $(d^2+(d+\varepsilon_2)^2/2)/(d+\varepsilon_2)=1.5d-0.5\varepsilon_2+O(\varepsilon_2^2)$, which implies that the average total cost of unit 3 is lower. Hence, in the solution for (16), unit 2 is offline, while unit 3 sets the price equal to its average total cost of output for a supply volume of $d+\varepsilon_3$. Since unit 3 has zero profit at this price, it can be either online or offline according to a solution for (16). As for the outputs of units 2 and 3 in the solution for the dual problem (6), we note that since the average total cost for output $g_2^{ec.\min}=\sqrt{2}d$ of unit 2 is above the average total cost for output $g_3^{ec.\min}=0.7d$ of unit 3, the solution of the dual problem (6) implies that unit 2 is offline, while unit 3 sets the market price and can have either status with zero profit.

**5. The set of market prices in the proposed approach**

If no LNMGU is present in the power system, then $P^+=\breve{P}(\varepsilon_{\hat{I}})$, $\forall \varepsilon_i>0$, $\forall i \in \hat{I}$. If at least one LNMGU is present in the power system, it has the following implications. If some LNMGU is online in the centralized dispatch, then the competitive equilibrium is impossible as the system has no equilibrium market price. For the CHP method, this implies nonzero duality gap and, consequently, nonzero total uplift payment. In the modified CHP approach, the online LNMGU needs no uplift payment if it is the only operating unit according to the centralized dispatch and receives nonzero uplift payment, otherwise. However, each LNMGU has zero profit in both CHP and the modified CHP methods. Also, each LNGMU $i$, $\forall i \in \hat{I}$, has the economic incentives to supply the output above demand if the market price is above $[w_i+c_i(g_i^{ec.\min})]/(g_i^{ec.\min})$ in CHP method (above $[w_i+c_i(d+\varepsilon_i)]/(d+\varepsilon_i)$ for a given $\varepsilon_i$ in the modified CHP), which provides an upper bound on the possible values of the market prices. Define $p^{LNMGU} = \min_{i \in \hat{I}} [w_i + c_i(g_i^{ec.\min})]/(g_i^{ec.\min})$, $\tilde{p}^{LNMGU}(\varepsilon_{\hat{I}}) = \min_{i \in \hat{I}} [w_i + c_i(d + \varepsilon_i)]/(d + \varepsilon_i)$. We have $p^+ \le p^{LNMGU}$, $\forall p^+ \in P^+$, and $\tilde{p}^+ \le \tilde{p}^{LNMGU}(\varepsilon_{\hat{I}})$, $\forall \tilde{p}^+ \in \breve{P}^+(\varepsilon_{\hat{I}})$. Thus, each of $P^+$ and $\breve{P}(\varepsilon_{\hat{I}})$ is bounded in the presence of at least one LNMGU and is either a singleton or a bounded closed interval. (Clearly, $p^{LNMGU} < \tilde{p}^{LNMGU}(\varepsilon_{\hat{I}})$, for sufficiently small $\forall \varepsilon_i > 0$, $\forall i \in \hat{I}$, which entails that the upper bound for $\breve{P}(\varepsilon_{\hat{I}})$ is higher than that for $P^+$.) Since for $\forall i \in \hat{I}$ both $\pi_i(p) = 0$, $\forall p \le p^{LNMGU}$, and $\tilde{\pi}_i(p, \varepsilon_i) = \hat{\pi}_i(p, \varepsilon_i) = 0$, $\forall p \le \tilde{p}^{LNMGU}(\varepsilon_{\hat{I}})$, we have

$$\nu^D = \max_{p \le p^{LNMGU}} [pd - \sum\nolimits_{i \in \breve{I}} \pi_i(p)], \; \tilde{\nu}^D(\varepsilon) = \hat{\nu}^D(\varepsilon) = \breve{\nu}^D(\varepsilon_{\hat{I}}) = \max_{p \le \tilde{p}^{LNMGU}(\varepsilon_{\hat{I}})} [pd - \sum\nolimits_{i \in \breve{I}} \pi_i(p)].$$

These two optimization problems are convex and have the identical (concave) objective functions. Since $p^{LNMGU} < \tilde{p}^{LNMGU}(\varepsilon_{\hat{I}})$, we have the following possibilities:

- Case 1: $p^{LNMGU} \notin P^+$. This implies that all maximum points of $pd - \sum\nolimits_{i \in \breve{I}} \pi_i(p)$ on $p \le \tilde{p}^{LNMGU}(\varepsilon_{\hat{I}})$ satisfy $p < p^{LNMGU}$. Thus, $P^+ = \breve{P}(\varepsilon_{\hat{I}})$, $\nu^D = \tilde{\nu}^D(\varepsilon)$, and both approaches give the same pricing outcome. In this case, all LNMGUs are irrelevant in both methods as the set of market prices is set by the other units.
- Case 2: $p^{LNMGU} \in P^+$, $P^+$ is not a singleton. This entails $P^+ \subset \breve{P}(\varepsilon_{\hat{I}})$ and $\nu^D = \tilde{\nu}^D(\varepsilon)$. Proposition 1 entails that $\nu = \nu^D$, and the both methods produce zero total uplift payments but the modified CHP may have a larger set of the market prices. We also have $\tilde{p}^+ \in P^+$ for $\forall \tilde{p}^+ \in \breve{P}(\varepsilon_{\hat{I}})$ such that $\tilde{p}^+ \le p^{LNMGU}$.
- Case 3: $P^+ = \{p^{LNMGU}\}$, $p^{LNMGU} \in \breve{P}(\varepsilon_{\hat{I}})$. Therefore, $P^+ \subset \breve{P}(\varepsilon_{\hat{I}})$ with $p^{LNMGU} \le \tilde{p}^+$, $\forall \tilde{p}^+ \in \breve{P}(\varepsilon_{\hat{I}})$. As a result, $\nu^D = \tilde{\nu}^D(\varepsilon)$ and the both methods produce the identical total uplift payments, but the modified CHP may have a larger set of the market prices.
- Case 4: $P^+ = \{p^{LNMGU}\}$, $p^{LNMGU} \notin \breve{P}(\varepsilon_{\hat{I}})$. In this case, $P^+ \cap \breve{P}(\varepsilon_{\hat{I}}) = \emptyset$, and we have $p^{LNMGU} < \tilde{p}^+$, $\forall \tilde{p}^+ \in \breve{P}(\varepsilon_{\hat{I}})$. Also, $\nu^D < \tilde{\nu}^D(\varepsilon)$, and the modified CHP approach produces the higher market price(s) and the lower total uplift payment.

Therefore, the modified CHP gives a lower total uplift payment only if $P^+ = \{p^{LNMGU}\}$, $p^{LNMGU} \notin \breve{P}(\varepsilon_{\hat{I}})$, and this reduction in the total uplift payment is accompanied by the market price increase. In the other cases, both methods produce the identical total uplift payments. Also, we have the following relation between the sets of market prices resulting from CHP and the proposed modified approach, which shows that the modified CHP tends to produce higher market prices.

*Proposition 6*: If $P^+ \cap \breve{P}(\varepsilon_{\hat{I}}) \neq \emptyset$, then $P^+ \subset \breve{P}(\varepsilon_{\hat{I}})$ and for $\forall \tilde{p}^+ \in \breve{P}(\varepsilon_{\hat{I}})$ the following holds for sufficiently small $\varepsilon_i > 0$, $\forall i \in \hat{I}$: either $\tilde{p}^+ \in P^+$ or $p^+ < \tilde{p}^+$, $\forall p^+ \in P^+$. If $P^+ \cap \breve{P}(\varepsilon_{\hat{I}}) = \emptyset$, then $p^+ < \tilde{p}^+$, $\forall p^+ \in P^+, \forall \tilde{p}^+ \in \breve{P}(\varepsilon_{\hat{I}})$.

*Proof*. If no LNMGU is present in the power system, then the proposition trivially holds. If at least one LNMGU is present, then the proposition follows from the above four cases. ■

Let us examine in more detail the structure of $\breve{P}(\varepsilon_{\hat{I}})$ for a power system that has at least one LNGMU. Define the reduced aggregate supply curve as the aggregate supply curve of all the generating units excluding $\hat{I}$, i.e. omitting all LNMGUs. The following three situations are possible.

In the first situation, the reduced aggregate supply curve has a point with the output volume lower than $d$ at a price $\min_{i \in \hat{I}} \left[ w_i + c_i(d) \right] / d$. For any sufficiently small $\varepsilon_i > 0$, $i \in \hat{I}$, we have $[w_i + c_i(d + \varepsilon_i)]/(d + \varepsilon_i) < [w_i + c_i(d)]/d$, hence any LNMGU from $\hat{I}_{\min}(\varepsilon_{\hat{I}})$ becomes marginal in the solution for (16) and sets the price. Thus, in this instance, $\breve{P}(\varepsilon_{\hat{I}})$ is a singleton with an element $\min_{i \in \hat{I}_{\min}(\varepsilon_{\hat{I}})} \left[ w_i + c_i(d + \varepsilon_i) \right] / (d + \varepsilon_i)$, which is a continuous function of $\varepsilon_i > 0$, $i \in \hat{I}$ as the minimum of a finite number of functions continuous in $\varepsilon_i$. Thus, as $\varepsilon_i \to +0$ for $\forall i \in \hat{I}$, the market price increases and attains the value of $\min_{i \in \hat{I}} \left[ w_i + c_i(d) \right] / d$. (This corresponds to the Case 4 above.) The second situation is realized when at a price $\min_{i \in \hat{I}} \left[ w_i + c_i(d) \right] / d$ the minimum supply volume on the reduced aggregate supply curve equals $d$. In this instance, we have two possibilities. If for any price below $\min_{i \in \hat{I}} \left[ w_i + c_i(d) \right] / d$ all the output volumes on the reduced aggregate supply curve are below $d$, then the conclusion for the first situation is applicable. Otherwise, for sufficiently small $\varepsilon_{\hat{I}}$, the set $\breve{P}(\varepsilon_{\hat{I}})$ is given by the bounded closed interval $[a, b(\varepsilon_{\hat{I}})]$ with $a$ being some nonnegative real number independent of $\varepsilon$ and $b(\varepsilon_{\hat{I}}) = \min_{i \in \hat{I}_{\min}(\varepsilon_{\hat{I}})} \left[ w_i + c_i(d + \varepsilon_i) \right] / (d + \varepsilon_i)$, which is continuous in $\varepsilon_{\hat{I}}$. As $\varepsilon_i \to +0$, $\forall i \in \hat{I}$, the set $\breve{P}(\varepsilon_{\hat{I}})$ tends to $[a, b(0)]$ with $b(0) = \min_{i \in \hat{I}} \left[ w_i + c_i(d) \right] / d$. (This matches the Cases 2, 3, 4 above.) Now we turn to the third situation, when at a price $\min_{i \in \hat{I}} \left[ w_i + c_i(d) \right] / d$, all the points on the reduced aggregate supply curve have output volumes higher than $d$. This means that for sufficiently small positive $\varepsilon_i$, $\forall i \in \hat{I}$, all LNMGUs are irrelevant for the dual problem (16) solution as the market price (a set of prices) in the modified CHP approach is fixed by the other units. (This corresponds to any of the Cases 1 – 4 above.) As a result, the set of prices in the modified CHP method has a well-defined limit as $\varepsilon \to +0$.

## 6. Examples

*Example 1 revisited.* Clearly, the generator is an LNMGU, and the generator modified private feasible set needs to be utilized in the proposed approach. Application of the modified CHP to Example 1 gives the singleton $\tilde{P}^+(\varepsilon)$ with the element $\tilde{p}^+(\varepsilon) = a + w/(d+\varepsilon)$ for sufficiently small $\varepsilon > 0$, which satisfies $0 < \varepsilon < g^{\max} - d$. Clearly, we have $\tilde{p}^+(\varepsilon) > p^+$. In this case, the generator receives the uplift payment equal $w[1 - d/(d+\varepsilon)]$, which is smaller than the uplift implied by CHP and is zero in the limit as $\varepsilon \to +0$. In both methods, the generator receives zero profit, but the market price resulting from the application of the proposed approach provides the proper economic signals for new entrants as it reflects the average cost of output needed to fully replace the incumbent generator. We note that contrary to the market price obtained from CHP, which is independent from $d$, the price $\tilde{p}^+(\varepsilon)$ is decreasing with the load.

*Example 3.* Let us consider a system with fixed demand $d$ and two generating units: unit 1 with $0 \le g_1 \le g_1^{\max}$ and the cost function $C_1(x_1) = w_1 u_1 + a_1 g_1$ and unit 2 with $0 \le g_2 \le g_2^{\max}$ and $C_2(x_2) = a_2 g_2$. The parameters are assumed to satisfy $a_1 > 0$, $a_2 > 0$, $a_1 + w_1 / g_1^{\max} < a_2 < a_1 + w_1 / d$, $d < g_1^{\max}$, $d < g_2^{\max}$. (Clearly, unit 1 is LNMGU while unit 2 is not.) These conditions ensure both that the primal problem (1) has a unique solution with unit 1 being offline (and unit 2 producing $d$) and that the dual problem (6) results in the unique market price $p^+ = a_1 + w_1 / g_1^{\max}$, which is set by unit 1. Thus, the market price is below $a_2$, which is the marginal cost of output of the online unit. The total uplift equals $(a_2 - p^+)d$ and is paid to unit 2 only. The application of the modified CHP for sufficiently small $\varepsilon_1 > 0$ results in the unique market price $\tilde{p}^+(\varepsilon_1, \varepsilon_2) = a_2$ (the set $\tilde{P}^+(\varepsilon_1, \varepsilon_2)$ is a singleton), which implies zero total uplift payment since unit 1 is not entitled for the uplift payment due to $\tilde{p}^+(\varepsilon_1, \varepsilon_2) < a_1 + w_1/(d + \varepsilon_1)$. Thus, $\tilde{p}^+(\varepsilon_1, \varepsilon_2) > p^+$ and all units have zero profits in both CHP and the proposed modified CHP. (We note that the proposed pricing method results in the market price being set by the online unit 2, not the offline unit 1. However, this is a specific property of the given example, and, in the general case, the new method also allows an offline unit to set the market price.)

*Example 4.* Let us amend Example 3 and replace the condition $d < g_2^{\max}$ by $d = g_2^{\max}$. Both the primal problem (1) solution for outputs and the dual problem (6) solution for the market price $p^+$ as well as each unit uplift payment in CHP remain the same. However, the outcome of the modified convex pricing method changes and is given by the bounded closed interval $\tilde{P}^+(\varepsilon_1, \varepsilon_2) = [a_2, a_1 + w_1/(d+\varepsilon_1)]$ for sufficiently small $\varepsilon_1 > 0$. We note that $p^+$ is below any elements of $\tilde{P}^+(\varepsilon_1, \varepsilon_2)$. In the limit as $\varepsilon_1 \to +0$, $\varepsilon_2 \to +0$, we have $\tilde{P}^+(+0,+0) = [a_2, a_1 + w_1/d]$ and zero total uplift. This manifests the fact that at a price above $a_1 + w_1/d$ unit 1 has an opportunity to sign profitable contracts with all the consumers to supply power volume $d$. As a result, if the market price were set above $a_1 + w_1/d$, then unit 1 would have to be compensated for the lost profit.

*Example 5.* Let us modify Example 3 and replace $d < g_2^{\max}$ by $d > g_2^{\max}$. Since unit 1 has to be online in the centralized dispatch problem and $a_1 < a_2$, we conclude that unit 2 is offline and unit 1 has output equal $d$ in the solution to (1). Due to $a_1 + w_1/g_1^{\max} < a_2$, CHP produces a unique market price $p^+ = a_1 + w_1/g_1^{\max}$, which implies the total uplift of $w(1 - d/g^{\max})$ paid to unit 1. The relation $a_2 < a_1 + w_1/d$ implies that we have $a_2 < a_1 + w_1/(d+\varepsilon_1)$ for sufficiently small $\varepsilon_1 > 0$, and the modified convex pricing algorithm results in the unique market price $\tilde{p}^+(\varepsilon_1, \varepsilon_2) = a_1 + w_1/(d+\varepsilon_1)$, which is higher than $p^+$. In the limit as $\varepsilon_1 \to +0$, $\varepsilon_2 \to +0$, no uplift is paid to unit 1, while unit 2 receives the uplift of $(a_1 + w_1/d - a_2)g_2^{\max}$. It is straightforward to verify that the modified CHP gives lower total uplift than the convex hull pricing method.

## 7. Discussion

In the presence of non-convexities, the equilibrium market price may not exist. In these cases, many features of the equilibrium market price become mutually exclusive, and any pricing method forfeits some of the equilibrium price attributes in favor of the other desired properties. In our approach, we stay within the uniform market pricing framework. To ensure the economic stability of the centralized dispatch outcome, we introduce both the uplift payments that compensate the lost profit calculated with respect to the opportunities attainable at the market (and the option not to participate in the market) and the financial penalties that provide the economic incentives not to explore the other opportunities. In the proposed pricing method, the opportunities available at the market include only those that correspond to a dispatch schedule (for all market players) satisfying both the system-level constraints and the market player private constraints as well as the economic constraints. The latter ensure that the resulting schedule cannot be optimized further by reducing the output of any market players. Thus, the proposed pricing scheme can be viewed as a mixture of the convex hull pricing and the nonlinear pricing approaches. We note that if the equilibrium price exists, then the proposed pricing algorithm gives the equilibrium price.

While CHP method results in the set of market prices that minimizes the total uplift payment in the uniform (linear) pricing framework, the major advantage of the proposed pricing scheme is further reduction of the total uplift payment. However, compared to the other pricing methods with a uniform market price, in our pricing algorithm the set of the external parameters, which are passed to a market player, includes not only the market price but also the additional constraints needed to define the opportunities available to the market player at the market. Our pricing approach also fits the general economic framework that each market player faces a rule (14) stating how its output revenue is calculated, and, given this rule, the market player lost profit is calculated based on the maximum value of its profit on the market player private feasible set.

The proposed approach is quite general and can be straightforwardly generalized to the multi-node multi-period power markets with price-sensitive bids and extended commodity space (such as power and reserve products). In these cases, the condition (9) with the proper objective function needs to include all the products and services traded at the market, the relevant private constraints of the market player (including the time-coupling constraints), and the system-level constraints (such as power flow and reserve adequacy constraints).

Strategic bidding analyses in both the convex hull pricing method and the proposed approach are complicated due to multi-part bids and non-convex private feasible sets of the market players. On one hand, the proposed method implies utilization of the reduced private feasible sets of some market players in the lost profit calculation, which may suppress the gaming strategies. On the other hand, the suggested pricing scheme may stimulate the improper behavior as a market player stated bid and technical parameters of its generating units may influence the additional constraints introduced in the private feasible sets of the other market participants and affect the resulting uplift payments as well as the market price. These questions have to be addressed in further research on the subject.

## 8. Conclusion

CHP treats each output volume allowed by generator private constraints as possible even if that output is technologically and/or economically infeasible. Thus, in this approach, all output volumes allowed by the generator private constraints enter the lost profit calculation. We propose modifying CHP by excluding from the lost profit computation all the nonzero output volumes that are technologically and/or economically infeasible. For this purpose, for each generator, we define a set of output volumes that are attainable at the centralized dispatch, i.e. both technologically and economically feasible in the centralized dispatch. These output volumes can be obtained as solutions to the centralized dispatch problem with the generator maximum output limits no higher than those in the original problem. For each generator, we define the modified private feasible set to include both its output volumes that are attainable as the centralized dispatch outcome and small $\varepsilon$- neighborhoods of these volumes within the generator original private feasible set (to indicate in the dual problem that the infinitesimal deviations of the market price from a given value may result in undersupply or oversupply of power) as well as the zero output volume (to ensure the non-confiscatory pricing). The procedure of reducing the generator private feasible set to the modified private feasible set can be viewed as extending the generator private constraint set to include new constraints depending on the status-output variable of this generator only. These constraints hold at the optimal point of the centralized dispatch problem and, therefore, do not affect the centralized dispatch outcome.

We explicitly construct the generator modified private feasible sets for a single-period uninode power system with fixed demand and zero generator minimum output. These sets are further used to calculate the set of the market prices and individual uplift payments in the limit of $\varepsilon \to +0$. To ensure the economic stability of the centralized dispatch outcome, in addition to the uplift payment for following the centralized dispatch we introduce the financial penalty for the output volumes outside the modified private feasible set of the generator.

The proposed modified CHP can be viewed as the producer revenue function amendment that includes both the uplift payment to a generator for following the centralized dispatch and the penalty for the output outside its modified private feasible set. The proposed approach produces a total uplift payment lower than (or equal to) the one resulting from CHP. The mathematical reason for potentially lower total uplift payment is the reduction of the market player private feasible sets used in the dual problem that produces both the set of market prices and the lost opportunity costs of market players. In the case of a convex centralized dispatch optimization problem, the proposed algorithm reproduces the marginal pricing.

Compared to CHP, only LNMGUs (the units with the maximum output limit above the fixed load and the decreasing average cost functions up to some output volume exceeding the fixed load) require special treatment in the decentralized dispatch problems emerging in the dual of the modified primal problem. Only for these units the modified private feasible sets need be used when calculating the set of market prices in the modified CHP. For the rest of the generating units, the original private feasible sets need no modifications, i.e. the decentralized dispatch problems for these units in the modified CHP method can be formulated the same way as in CHP. Consequently, if a power system has no LNMGU, then the proposed method gives the same set of market prices as CHP. (We note that this result heavily relies on the assumptions of zero generator minimum output limits and absence of the price-sensitive load in the power system.)

We showed that a set of market prices resulting from the proposed pricing method has a well-defined limit as $\varepsilon \to +0$. In addition, the new approach gives a set of market prices that tends to be "no lower" than the set of prices obtained from CHP.

## Appendix

*Proof of Proposition 1*: Let $P^+$ contain more than one element. Since $P^+$ is a convex set, there is some $p' \in \text{int } P^+$. Clearly, $(p'-\alpha, p'+\alpha) \in \text{int } P^+$ for some $\alpha > 0$. From (6) it follows that on $(p'-\alpha, p'+\alpha)$ the convex function $S(p) = \sum_{i \in I} \pi_i(p) - pd$ is constant with $\partial S(p) = \{0\}$. From $\partial S(p) = Conv\{\sum_{i \in I} g_i - d \mid \nu^D = p(d - \sum_{i \in I} g_i) + \sum_{i \in I} C_i(x_i), x_i \in X_i, \forall i \in I\}$, we conclude that $\forall p \in (p'-\alpha, p'+\alpha) \subset P^+$ the power balance constraint $\sum_{i \in I} g_i = d$ holds for $\forall x \in \arg \min_{x_i \in X_i, \forall i \in I} [p(d - \sum_{i \in I} g_i) + \sum_{i \in I} C_i(x_i)]$. As a results, $\nu = \nu^D$. ■

*Proof of Proposition 3*: If for a given $i$ we have $\overline{g}_i^{\min} \le \varepsilon_i$, then $\tilde{X}_i(\varepsilon_i) = \hat{X}_i(\varepsilon_i)$ and, hence, $\tilde{\pi}_i(p, \varepsilon_i) = \hat{\pi}_i(p, \varepsilon_i)$, $\forall p \in R$. Therefore, in this case, $\tilde{\pi}_i(p, \varepsilon_i)$ can be replaced by $\hat{\pi}_i(p, \varepsilon_i)$ in (12) with no effect on both $\tilde{P}^+(\varepsilon)$ and $\tilde{\nu}^D(\varepsilon)$. Thus, we focus on $i$ with $\overline{g}_i^{\min} > \varepsilon_i$. Let us define for $g_i \in R$ an extended real-valued function $\tilde{f}_i(g_i)$ as follows: $\tilde{f}_i(g_i) = f_i(g_i)$ for $g_i$ that corresponds to the output volumes in $\tilde{X}_i(\varepsilon_i)$ and $\tilde{f}_i(g_i) = +\infty$, otherwise; likewise, we define an extended real-valued function $\hat{f}_i(g_i)$ using output volumes in $\hat{X}_i(\varepsilon_i)$. For $g_i \le 0$ and $g_i \ge \overline{g}_i^{\min} - \varepsilon_i$, we have $\tilde{f}_i(g_i) = \hat{f}_i(g_i)$. Consider their respective convex hulls $\tilde{f}_i^h(g_i)$ and $\hat{f}_i^h(g_i)$, which have finite values on $0 \le g_i \le \min(\overline{g}_i^{\max} + \varepsilon_i; g_i^{\max})$. It is straightforward to verify that for $g_i \ge \overline{g}_i^{\min} - \varepsilon_i$ we have $\tilde{f}_i^h(g_i) = \hat{f}_i^h(g_i)$. Thus, $\partial\tilde{f}_i^h(g_i) = \partial\hat{f}_i^h(g_i)$ for $g_i > \overline{g}_i^{\min} - \varepsilon_i$. Clearly, the function $\hat{\pi}_i(p, \varepsilon_i)$ is convex in $p$ (as a point-wise maximum of the function linear in $p$) with the subdifferential with respect to $p$ given by

$$\partial\hat{\pi}_i(p, \varepsilon_i) = \{g_i \mid g_i \in dom \quad \hat{f}_i^h(g_i), p \in \partial\hat{f}_i^h(g_i)\}. \quad (17)$$

The analogous expression holds for $\partial\tilde{\pi}_i(p, \varepsilon_i)$. Consequently, we conclude that $\forall p \in R$ the sets $\partial\tilde{\pi}_i(p, \varepsilon_i)$ and $\partial\hat{\pi}_i(p, \varepsilon_i)$ when restricted to $(\overline{g}_i^{\min} - \varepsilon_i, \min(\overline{g}_i^{\max} + \varepsilon_i; g_i^{\max})]$ with $\varepsilon_i > 0$ are identical (possibly empty, for some $p \in R$). As (15) is a maximization problem of a concave function, its set of maximizers $\hat{P}^+(\varepsilon)$ is given by $d \in \sum_{i \in I} \partial\hat{\pi}_i(p, \varepsilon_i)$. The definition of $\overline{g}_i^{\min}$ implies that the elements of $\partial\hat{\pi}_i(p, \varepsilon_i)$ lower than $\overline{g}_i^{\min}$ (if any) do not affect the set $\hat{P}^+(\varepsilon)$. (Also, if all elements of $\partial\hat{\pi}_i(p, \varepsilon_i)$ are below $\overline{g}_i^{\min}$, then $d \notin \sum_{i \in I} \partial\tilde{\pi}_i(p, \varepsilon_i)$ and $p \notin \tilde{P}^+(\varepsilon)$.) An analogous observation holds for $\tilde{P}^{+}(\varepsilon)$. Thus, when $\varepsilon_i > 0$, $\forall i \in I$, the solutions to $d \in \sum_{i \in I} \partial \tilde{\pi}_i(p, \varepsilon_i)$ and $d \in \sum_{i \in I} \partial \hat{\pi}_i(p, \varepsilon_i)$ are identical, giving us $\hat{P}^{+}(\varepsilon) = \tilde{P}^{+}(\varepsilon)$. Since both $pg_i - \tilde{f}_i^{h}(g_i)$ and $pg_i - \hat{f}_i^{h}(g_i)$ are the concave functions on the convex set $0 \le g_i \le \min(\bar{g}_i^{\max} + \varepsilon_i; g_i^{\max})$, the sets $\partial \tilde{\pi}_i(p, \varepsilon_i)$ and $\partial \hat{\pi}_i(p, \varepsilon_i)$ are the sets of maximizers of $\tilde{\pi}_i(p, \varepsilon_i)$ and $\hat{\pi}_i(p, \varepsilon_i)$, respectively. Since $\tilde{P}^{+}(\varepsilon)$ is nonempty, the restrictions on $[\bar{g}_i^{\min}, \min(\bar{g}_i^{\max} + \varepsilon_i; g_i^{\max})]$ of the both sets $\partial \tilde{\pi}_i(p, \varepsilon_i)$ and $\partial \hat{\pi}_i(p, \varepsilon_i)$ are nonempty, $\forall p \in \tilde{P}^{+}(\varepsilon)$. The identity of $\partial \tilde{\pi}_i(p, \varepsilon_i)$ and $\partial \hat{\pi}_i(p, \varepsilon_i)$ on $(\bar{g}_i^{\min} - \varepsilon_i, \min(\bar{g}_i^{\max} + \varepsilon_i; g_i^{\max})]$ along with $\tilde{f}_i^{h}(g_i) = \hat{f}_i^{h}(g_i)$ for $\bar{g}_i^{\min} \le g_i < +\infty$ entails $\hat{\pi}_i(p, \varepsilon_i) = \tilde{\pi}_i(p, \varepsilon_i)$, $\forall p \in \tilde{P}^{+}(\varepsilon)$, for the case $\bar{g}_i^{\min} > \varepsilon_i$. The statements of the first two bullets trivially imply the claim of the third bullet.■

*Proof of Proposition 4*: Clearly, if for a given $i \in \breve{I}$ we have $\bar{g}_i^{\max} + \varepsilon_i \ge g_i^{\max}$, then $X_i = \hat{X}_i(\varepsilon_i)$, $\pi_i(p) = \hat{\pi}_i(p, \varepsilon_i)$, $\forall p \in R$. In this case, $\hat{\pi}_i(p, \varepsilon_i)$ can be replaced by $\pi_i(p)$ in (15) with no effect on $\hat{P}^{+}(\varepsilon)$ or $\hat{v}^{D}(\varepsilon)$. Thus, we may restrict our consideration to the case of $\bar{g}_i^{\max} + \varepsilon_i < g_i^{\max}$, $i \in \breve{I}$. It is straightforward to verify that in this case $f_i^{h}(g_i) = \hat{f}_i^{h}(g_i)$ for $g_i \le \bar{g}_i^{\max} + \varepsilon_i$ and, hence, $\partial f_i^{h}(g_i) = \partial \hat{f}_i^{h}(g_i)$ for $g_i < \bar{g}_i^{\max} + \varepsilon_i$. Using (4) and (17) for $\forall i \in \breve{I}$, we conclude that when restricted to $[0, \bar{g}_i^{\max} + \varepsilon_i)$, the sets $\partial \pi_i(p)$ and $\partial \hat{\pi}_i(p, \varepsilon_i)$ are identical, $\forall p \in R$. Since (16) is a maximization problem of a concave objective function, $\breve{P}^{+}(\varepsilon_{\hat{I}})$ is given by the set of solutions to $d \in \sum_{i \in \breve{I}} \partial \pi_i(p) + \sum_{i \in \hat{I}} \partial \hat{\pi}_i(p, \varepsilon_i)$, and, as a result, the elements of $\partial \pi_i(p)$ higher than $\bar{g}_i^{\max}$ (if any) do not affect the set $\breve{P}^{+}(\varepsilon_{\hat{I}})$. Therefore, for $\varepsilon_i > 0$, $\forall i \in \breve{I}$, the equations $d \in \sum_{i \in \breve{I}} \partial \pi_i(p) + \sum_{i \in \hat{I}} \partial \hat{\pi}_i(p, \varepsilon_i)$ and $d \in \sum_{i \in I} \partial \hat{\pi}_i(p, \varepsilon_i)$ have identical sets of solutions for $p$, which entails $\breve{P}^{+}(\varepsilon_{\hat{I}}) = \hat{P}^{+}(\varepsilon)$. Since the elements of $\partial \hat{\pi}_i(p, \varepsilon_i)$ higher than $\bar{g}_i^{\max}$ (if any) do not contribute to $d \in \sum_{i \in I} \partial \hat{\pi}_i(p, \varepsilon_i)$, (we still consider the case of $\bar{g}_i^{\max} + \varepsilon_i < g_i^{\max}$, $i \in \breve{I}$), we conclude that $\forall p \in \hat{P}^{+}(\varepsilon)$ the set $\partial \hat{\pi}_i(p, \varepsilon_i)$ has at least one element in $[0, \bar{g}_i^{\max}]$; hence, $\forall p \in \hat{P}^{+}(\varepsilon)$, both $\partial \hat{\pi}_i(p, \varepsilon_i)$ and $\partial \pi_i(p)$ are nonempty and equal when restricted to this range. As both $pg_i - f_i^{h}(g_i)$ and $pg_i - \hat{f}_i^{h}(g_i)$ are the concave functions on the convex set $0 \le g_i \le \bar{g}_i^{\max} + \varepsilon_i$, the sets $\partial \pi_i(p)$ and $\partial\hat{\pi}_i(p,\varepsilon_i)$ are the maximizers of $\pi_i(p)$ and $\hat{\pi}_i(p,\varepsilon_i)$, respectively. Since $f_i^h(g_i)=\hat{f}_i^h(g_i)$ for $0\le g_i\le \bar{g}_i^{\max}+\varepsilon_i$, we obtain $\hat{\pi}_i(p,\varepsilon_i)=\pi_i(p)$, $\forall i\in\breve{I}$, $\forall p\in\hat{P}^+(\varepsilon)$, which also entails $\breve{v}^D(\varepsilon_{\hat{I}})=\hat{v}^D(\varepsilon)$. ■

*Proof of Proposition 5*: We have

$$d\in\sum\nolimits_{i\in\breve{I}}\partial\pi_i(p)+\sum\nolimits_{i\in\hat{I}}\partial\hat{\pi}_i(p,\varepsilon_i)\,,\ \forall p\in\breve{P}^+(\varepsilon_{\hat{I}})\,. \quad (18)$$

For sufficiently small $\varepsilon_i>0$, $i\in\hat{I}$, the application of (4) yields:

$$\partial\hat{\pi}_i(p,\varepsilon_i)=\begin{cases}0,\ p<[w_i+c_i(d+\varepsilon_i)]/(d+\varepsilon_i)\\ [0,d+\varepsilon_i],\ p=[w_i+c_i(d+\varepsilon_i)]/(d+\varepsilon_i)\\ d+\varepsilon_i,\ p>[w_i+c_i(d+\varepsilon_i)]/(d+\varepsilon_i)\end{cases}. \quad (19)$$

Hence, if for some $i\in\hat{I}$ we have $p>[w_i+c_i(d+\varepsilon_i)]/(d+\varepsilon_i)$, then $\partial\hat{\pi}_i(p,\varepsilon_i)=d+\varepsilon_i>d$, and since all elements of $\partial\pi_i(p)$ are nonnegative, (18) implies $p\notin\breve{P}^+(\varepsilon_{\hat{I}})$. Therefore, $p\le\min_{i\in\hat{I}}\,[w_i+c_i(d+\varepsilon_i)]/(d+\varepsilon_i)$, $\forall p\in\breve{P}^+(\varepsilon_{\hat{I}})$, and all LNMGUs with $[w_i+c_i(d+\varepsilon_i)]/(d+\varepsilon_i)$ above $\min_{i\in\hat{I}}\,[w_i+c_i(d+\varepsilon_i)]/(d+\varepsilon_i)$ do not contribute to RHS of (18) and can be disregarded in (16) without affecting $\breve{P}^+(\varepsilon_{\hat{I}})$. Also, if there is more than one LNMGU $i$ with $i\in\hat{I}_{\min}(\varepsilon_{\hat{I}})$, then all but any one of such LNMGUs can be removed from the consideration in (16) with no effect on $\breve{P}^+(\varepsilon_{\hat{I}})$. ■

## Declarations

- The author has no relevant financial or non-financial interests to disclose.
- The author has no conflicts of interest to declare that are relevant to the content of this article.
- The author received no funding to support this research.
- The author confirms that the data supporting the findings of this study are available within the article.